\documentstyle[12pt,leqno]{amsart}
\newtheorem{dummy}{}[section]
\newtheorem{definition}[dummy]{Definition}
\newtheorem{theorem}[dummy]{Theorem}

\newtheorem{proposition}[dummy]{Proposition}
\newtheorem{lemma}[dummy]{Lemma}
\newtheorem{corollary}[dummy]{Corollary}
\newtheorem{remark}[dummy]{Remark}

\begin{document}
\bibliographystyle{plain}
\title{Derived Langlands IV: Notes on ${\mathcal M}_{c}(G)$-induced representations}
\author{Victor P. Snaith}
\date{July 2020}
\maketitle
\tableofcontents

\section{Extending the definition of admissibility}

If $G$ is a locally profinite group and $k$ is an algebraically closed field then a $k$-representation of $G$ is a vector space $V$ with a left, $k$-linear $G$-action. Let ${\mathcal M}_{c}(G)$ denote the poset of 
pairs $(H, \phi)$ where $H$ is a subgroup of $G$, which is compact open and $\phi : H \longrightarrow k^{*}$ is a continuous character.

A representation $V$ is called smooth (\cite{BH06} p.13) if
\[ V = \bigcup_{K \subset G,  \ K \ {\rm compact, open}} \ V^{K} .\]
$V$ is called admissible (\cite{BH06} p.13) if ${\dim}_{k}( V^{K}) < \infty$ for all compact open subgroups $K$. 

Notice that, if $v_{i} \in V^{K_{i}}$ for $i=1,2$ with $K_{i}$ compact open then
a linear combination $a_{1} v_{1} + a_{2} v_{2} \in V^{K_{1} \bigcap K_{2}} $ for $a_{i} \in k$ with 
$K_{1} \bigcap K_{2}$ compact open so that 
\[ V =  \bigcup_{K \subset G,  \ K \ {\rm compact, open}} \ V^{K}  = {\rm Span}_{K \subset G,  \ K \ {\rm compact, open}} \  V^{K}.  \]

The smooth representations of $G$ form an abelian category.

Define a subspace of $V$, denoted by $V^{(H, \phi)}$, for $(H, \phi) \in {\mathcal M}_{c}(G)$ by
 \[  V^{(H, \phi)} = \{ v \in V \ | \ g \cdot v =  \phi(g)v \ {\rm for \ all} \ g \in H \}   . \]
 Hence $V^{K} = V^{( K,  1)}$ if $1$ denotes the trivial character.
 
 We shall say that $V$ is ${\mathcal M}_{c}(G)$-smooth\footnote{In \cite{Sn20} I made a childish typographical slip-up, writing $\bigcup$ instead of $\Sigma$ in this definition. Hence this time I have written ${\rm Span}$. In my current time-shortage (\cite{JMOS18} and \cite{XX})  the reader will probably find me making a number of (I hope inconsequential) slip ups.}  if 
 \[ V = {\rm Span}_{ (H, \phi) \in {\mathcal M}_{c}(G)} \ V^{(H, \phi)}.\]
 In addition we shall say that $V$ is ${\mathcal M}_{c}(G)$-admissible if 
 ${\rm dim}_{k}V^{(H, \phi)} < \infty$ for all  $(H, \phi) \in {\mathcal M}_{c}(G)$.

 \begin{proposition}{$_{}$}
 \label{1.1}
 \begin{em}
 
 Let $G$ be a locally profinite group and let $k$ be an algebraically closed field. Let $V$ be a $k$-representation of $G$. Suppose that every continuous, $k$-valued character of a compact open subgroup of $G$ has finite image. Then $V$ is 
 ${\mathcal M}_{c}(G)$-admissible if and only if it is admissible.
 \end{em}
 \end{proposition}
 
 Let us begin by recalling, from (\cite{Sn18} Chapter Two \S1), induced and compactly induced smooth representations.

\begin{definition}{Smooth induction}
\label{1.2}
\begin{em}

Let $G$ be a locally profinite group and $H \subseteq G$ a closed subgroup. Thus $H$ is also  locally profinite. Let 
\[   \sigma : H \longrightarrow  {\rm Aut}_{k}(W)    \]
be a smooth representation of $H$. Set $X$ equal to the space of functions $f: G \longrightarrow W  $ such that (writing simply $h \cdot w$ for $\sigma(h)(w)$ if $h \in H, w \in W$)

(i) \  $f(hg) = h \cdot f(g)$ for all $h \in H, g \in G$,

(ii)  \  there is a compact open subgroup $K_{f} \subseteq G$ such that $f(gk) = f(g)$ for all $g \in G, k \in K_{f}$.

The (left) action of $G$ on $X$ is given by $(g \cdot f)(x)= f(xg)$ and
\[   \Sigma :  G  \longrightarrow  {\rm Aut}_{k}(X)  \]
gives a smooth representation of $G$.

The representation $\Sigma$ is called the representation of $G$ smoothly induced from $\sigma$ and is usually denoted by $\Sigma = {\rm Ind}_{H}^{G}(\sigma)$.
\end{em}
\end{definition}

\begin{definition}{${\mathcal M}_{c}(G)$-smooth induction}
\label{1.3}
\begin{em}

Let $G$ be a locally profinite group and $H \subseteq G$ a closed subgroup. Let 
\[   \sigma : H \longrightarrow  {\rm Aut}_{k}(W)    \]
be an ${\mathcal M}_{c}(H)$-smooth representation of $H$. Set $X$ equal to the space of functions $f: G \longrightarrow W  $ such that (writing simply $h \cdot w$ for $\sigma(h)(w)$ if $h \in H, w \in W$)

(i) \  $f(hg) = h \cdot f(g)$ for all $h \in H, g \in G$,

(ii)  \  Each $f \in X$ may be written as a linear combination of the form $f = \sum_{i=1}^{n} \  a_{i}f_{i}$ where $a_{i} \in k$ and for each $1 \leq i \leq n$ 
there is $(K_{f_{i}}, \phi_{f_{i}}) \in {\mathcal M}_{c}(G)$ such that 
$f_{i}(gk) = \phi_{f_{i} }(k) f_{i} (g)$ for all $g \in G, k \in K_{f_{i}}$.

The (left) action of $G$ on $X$ is given by $(g \cdot f)(x)= f(xg)$ and
\[   \Sigma_{{\mathcal M}_{c}} :  G  \longrightarrow  {\rm Aut}_{k}(X)  \]
gives an ${\mathcal M}_{c}(G)$-smooth representation of $G$.

The representation $\Sigma_{{\mathcal M}_{c}} $ is called the representation of $G$ ${\mathcal M}_{c}$-smoothly induced from $\sigma$ and is will be denoted by $\Sigma_{{\mathcal M}_{c}} = {\rm IND}_{ H}^{G}(\sigma)$. 
\end{em}
\end{definition}

We just pause to check the following:
\begin{lemma}{$_{}$}
\label{1.4}
\begin{em}
 If $f \in X$ in Definition \ref{1.3} and $g \in G$ then $g \cdot f \in X$ also.
\end{em}
\end{lemma}  

{\bf Proof:} 

Suppose that $f \in X$ and $(K_{f_{i}}, \phi_{f_{i}}) \in {\mathcal M}_{c}(G)$ satisfy the second condition of Definition \ref{1.3}. Define $g( (K_{f_{i}}, \phi_{f_{i}}))$ to be the pair consisting of the group $gK_{f_{i}}g^{-1}$ with character 
 \[   (g^{-1})^{*}(\phi_{f_{i}}): gK_{f_{i}}g^{-1} \longrightarrow  k^{*} \]
where, for $k \in K_{f_{i}}$,
we set $(g^{-1})^{*}(\phi_{f_{i}})(gkg^{-1}) = \phi_{f_{i}}(k)$. Therefore 
\[  \begin{array}{ll}
(g \cdot f_{i})(g'gkg^{-1}) &  = f_{i}(g'gkg^{-1}g) \\
\\
& = f_{i}(g'gk) \\
\\
& = \phi_{f_{i}}(k) f_{i}(g'g) \\
\\
& = (g^{-1})^{*}(\phi_{f_{i}})(gkg^{-1}) f_{i}(g'g) \\
\\
& =  (g^{-1})^{*}(\phi_{f_{i}})(gkg^{-1})   (g \cdot f_{i})(g')  ,
\end{array}   \]
as required. Clearly action by $g \in G$ preserves condition (i) of Definition \ref{1.3} $\Box$

\begin{dummy}{The c-IND variation}
\label{1.5}
\begin{em}

Inside $X$ let $X_{c}$ denote the set of functions which are compactly supported modulo $H$. This means that the image of the support
\[   {\rm supp}(f) = \{ g \in G \ | \  f(g) \not= 0  \}  \]
has compact image in $H \backslash G$. Alternatively  there is a compact subset $C \subseteq G$ such that $  {\rm supp}(f) \subseteq H \cdot C$.

The $\Sigma_{{\mathcal M}_{c}}$-action on $X$ preserves $X_{c}$, since $ {\rm supp}(g \cdot f) =  {\rm supp}(f) g^{-1} \subseteq HCg^{-1}$, and we obtain $X_{c} =  c- {\rm IND}_{H}^{G}(W) $, the compact induction of $W$ from $H$ to $G$ \footnote{I have not yet bothered to check the following:  Note that this condition requires, in Definition \ref{1.3}, at first sight, that $f$ is compactly supported modulo $H$ but not necessarily the individual $f_{i}$'s.}.
\end{em}
\end{dummy}

 \begin{lemma}{$_{}$}
 \label{1.6}
 \begin{em}
 
 (i) \ If $K$ is an open subgroup of $G$ and $V$ is an  ${\mathcal M}_{c}(G)$-smooth representation then ${\rm Res}_{H}^{G}(V)$ is an  ${\mathcal M}_{c}(H)$-smooth representation.
 
  (ii) \  If $\chi : G \longrightarrow k^{*}$ is a continuous character and $V$ is an 
  ${\mathcal M}_{c}(G)$-admissible representation then $V \otimes_{k} k_{\chi}$ is also an ${\mathcal M}_{c}(G)$-admissible representation. Herre $k_{\chi}$ is $k$ acted upon by $G$ via the character $\chi$.
 \end{em}
 \end{lemma}
 
 {\bf Proof:}
 
 (i) \  If $v \in V$ lies in $V^{(J, \phi)}$ for some $J$ compact open in $G$ then $v \in V^{(H \bigcap J, \phi)}$ and $H \bigcap J$ is compact open in $H$. Therefore $V$ is ${\mathcal M}_{c}(H)$-smooth.
 
 (ii) \ This follows from the relation $(V \otimes_{k} k_{\chi})^{(H, \phi)} = V^{(H, \phi \otimes \chi^{-1})}$.
 
\begin{remark}{$_{}$}
\label{1.7}
\begin{em}
 Admissibility is not inherited by subgroups. For example, if $H$ is a compact open subgroup of the centre of $G$ and if $V$ has a central character $\underline{\phi}$ then ${\rm Res}_{H}^{G}(V)^{(H, \underline{\phi})} = V$.
\end{em}
\end{remark}
\begin{lemma}{$_{}$}
\label{1.8}
\begin{em}

Let 
\[ 0 \longrightarrow V_{1} \longrightarrow V_{2} \longrightarrow V_{3} \longrightarrow 0 \]
be a short exact sequence of $G$-representations.

If $K$ is compact open and $(K , \phi) \in {\mathcal M}_{c}(G)$ then
\[ 0 \longrightarrow V_{1}^{(K , \phi)} \longrightarrow V_{2}^{(K , \phi)} \longrightarrow V_{3}^{(K , \phi)} \longrightarrow 0 \]
is exact.
\end{em}
\end{lemma}

{\bf Proof:}

\[ 0 \longrightarrow V_{1} \otimes_{k} k_{\phi^{-1}} \longrightarrow V_{2} \otimes_{k} k_{\phi^{-1}}  \longrightarrow 
V_{3} \otimes_{k} k_{\phi^{-1}}  \longrightarrow 0 \]
is an exact sequence of $K$-representations. Applying $K$-fixed points to this sequence is exact, because $K$ is compact, and the result is the short exact sequence of $k$-vector spaces in the statement of the lemma. $\Box$

\begin{corollary}{$_{}$}
\label{1.9}
\begin{em}

Let 
\[ 0 \longrightarrow V_{1} \longrightarrow V_{2} \longrightarrow V_{3} \longrightarrow 0 \]
be a short exact sequence of $G$-representations. 

Then $V_{2}$ is 
${\mathcal M}_{c}(G)$-admissible if and only if $V_{1}$ and $V_{3}$ are both 
\linebreak
${\mathcal M}_{c}(G)$-admissible.
\end{em}
\end{corollary}
\begin{proposition}{$_{}$}
\label{1.10}
\begin{em}

Let $V$ be a smooth representation of $G$. Suppose that $K$ is a compact open subgroups and that
 $(K, \phi) \in {\mathcal M}_{c}(G)$. Then ${\rm dim}(V^{( K, \phi) }) < \infty$.
\end{em}
\end{proposition} 
{\bf Proof:}

By (\cite{WCnotes} Proposition 2.1.4 p.20), the restriction of $V$ to $K$ is the direct sum of finite-dimensional irreducible $K$-representations, each appearing with finite multiplicity. If this $K$-representation is $\oplus_{\alpha} \ V_{\alpha}$ then $V^{(K, \phi) } = \oplus_{\alpha, \ V_{\alpha} = k_{\phi}}  \ V_{\alpha}$ which is finite-dimensional. $\Box$

\begin{theorem}{$_{}$}
\label{1.11} 
\begin{em}

Let 
\[ 0  \longrightarrow  V_{1} \stackrel{i}{\longrightarrow} V_{2} \stackrel{j}{\longrightarrow} V_{3} \longrightarrow 0 \]
be a short exact sequence of $G$-representations. Then
$V_{2}$ is ${\mathcal M}_{c}(G)$-smooth if and only if $V_{1}$ and $V_{3}$ are both ${\mathcal M}_{c}(G)$-smooth.
\end{em}
\end{theorem}

{\bf Proof:}

Assume that $V_{2}$ is ${\mathcal M}_{c}(G)$-smooth. Therefore 
\[  V_{2} = {\rm Span}_{ (H, \phi) \in {\mathcal M}_{c}(G)} \ V_{2}^{(H, \phi)}.\]
Since $j( V_{2}^{(H, \phi)} ) \subseteq V_{3}^{(H, \phi)}$ and $j$ is surjective we see that 
\[  V_{3} = {\rm Span}_{ (H, \phi) \in {\mathcal M}_{c}(G)} \ V_{3}^{(H, \phi)}.\]

Now suppose that $0 \not= w \in V_{1}$ satisfies $w = u_{1} + u_{2} + \ldots + u_{n}$ with 
$u_{i} \in V_{2}^{(H_{i}, \phi_{i})}$ with $(H_{i}, \phi) \in {\mathcal M}_{c}(G)$. We may assume that $H = H_{1} = \ldots = H_{n}$ by taking intersections. We shall proceed by induction on $n$ and start the induction by 
Lemma \ref{1.8}. By induction we may assume that $w = u_{1} + u_{2} + \ldots + u_{n}$ with 
$u_{i} \in V_{2}^{(H, \phi_{i})}$ with $(H, \phi) \in {\mathcal M}_{c}(G)$ and the $\phi_{i}$'s are pairwise distinct.
We have $j(u_{1} + u_{2} + \ldots + u_{n}) = j(i(w)) = 0$. Choose $g \in G$ such that $\phi_{1}(g) \not= \phi_{2}(g)$. Then
\[  i(g \cdot w) = \phi_{1}(g) u_{1} + \phi_{2}(g) u_{2} + \ldots + \phi_{n}(g) u_{n} \]
so that 
\[  i(g \cdot w - \phi_{1}(g) w ) =  (\phi_{2}(g) - \phi_{1}(g) u_{2} + \ldots + (\phi_{n}(g) - \phi_{1}(g) u_{n} . \]
By induction, since 
 \[ j(  (\phi_{2}(g) - \phi_{1}(g) u_{2} + \ldots + (\phi_{n}(g) - \phi_{1}(g) u_{n}  ) = 0 ,\]
 we have 
 \[  g \cdot w - \phi_{1}(g) w  \in {\rm Span}_{2 \leq s \leq n} \ V_{1}^{(H, \phi_{s})} .\]
 Similarly 
  \[  g \cdot w - \phi_{2}(g) w  \in {\rm Span}_{s=1 \ {\rm or} \ 3 \leq s \leq n} \ V_{1}^{(H, \phi_{s})} .\]
 Therefore
 \[  \phi_{2}(g) w - \phi_{1}(g) w  \in  {\rm Span}_{1 \leq s \leq n} \ V_{1}^{(H, \phi_{s})}  \]
 and, since $ \phi_{2}(g)  \not=  \phi_{1}(g)$ we see that
  \[   w  \in  {\rm Span}_{1 \leq s \leq n} \ V_{1}^{(H, \phi_{s})},  \]
  as required. Therefore $V_{1}$ and $V_{3}$ are ${\mathcal M}_{c}(G)$-smooth.
 
 The converse follows from Lemma \ref{1.8}. $\Box$
 
 \begin{corollary}{$_{}$}
 \label{1.12}
 \begin{em}
 
 ${\mathcal M}_{c}(G)$-smooth representations form an abelian category of which ${\mathcal M}_{c}(G)$-admissible representations form an abelian subcategory.
 \end{em}
 \end{corollary}
 
 \begin{theorem}{$_{}$}
 \label{1.13}
 \begin{em}
  Let $H$ be a closed subgroup of $G$ and let $(\tau, U)$ be an ${\mathcal M}_{c}(H)$-smooth $H$-representation. Then
  
  (i) \ ${\rm IND}_{H}^{G}(\tau)$ is ${\mathcal M}_{c}(G)$-smooth.
  
  (ii) \ The map $\Lambda :  {\rm IND}_{H}^{G}(\tau)  \longrightarrow  U$ defined by $\Lambda(f) = f(1)$ is a surjective $H$-map.
  
  (iii) \ The restriction of $\Lambda$ to any non-zero $G$-subspace of ${\rm IND}_{H}^{G}(\tau)$ is non-trivial.
  
  (iv) \ If $H \backslash G$ is compact and $(\tau, U)$ is ${\mathcal M}_{c}(H)$-admissible then ${\rm IND}_{H}^{G}(\tau)$ is ${\mathcal M}_{c}(G)$-admissible.
  
  (v) \ If $(\pi, V)$ is ${\mathcal M}_{c}(G)$-smooth then composition with $\Lambda$ induces an isomorphism of the form
  \[  ( \Lambda \cdot - ) : {\rm Hom}_{G}(V, IND_{H}^{G}(\tau)) \stackrel{\cong}{\longrightarrow }  {\rm Hom}_{H}(V, U).  \]
 \end{em}
 \end{theorem}
 
 {\bf Proof:}
 
Part (i) is immediate from the definition.

For part (ii) we observe that $\Lambda( h \cdot f) = f(h) = \tau(h) f(1) = \tau(h)(\Lambda(f))$. Given $u \in U$ 
choose a compact open subgroups $K_{i}$ of $G$ such that 
\linebreak
$(K \bigcap H, \phi_{i}) \in {\mathcal M}_{c}(H)$ and $ u = \sum_{i} \ u_{i}$ with $u_{i} \in U^{(K \bigcap H, \phi_{i})}$. Define $f_{i} \in IND_{H}^{G}(\tau)^{(K \bigcap H, \phi_{i})}$ by
$f_{i}(g) =  \tau(h) u_{i}) $ if $ g=hk, h \in H, k \in K$ and  $f_{i}(g)=0$ otherwise. This is compactly supported modulo $H$ and $\Lambda(f_{i}) = u_{i}$ so $u = \Lambda( \sum_{i} f_{i})$, as required.

For part (iii) suppose that $V$ is a non-trivial $G$-subspace of ${\rm IND}_{H}^{G}(\tau)$ and choose $f \in V$ such that $f(g) \not= 0$ for some $g \in G$. Therefore $g \cdot f \in V$ also and $\Lambda( g \cdot f) = f(g)$.
 
 For part (iv) suppose that $K$ is a compact open subgroup of $G$ such that $(K, \phi) \in {\mathcal M}_{c}(G)$. Suppose that $X$ is a finite subset of $G$ and that $U_{0}$ is a finite dimensional subspace of $U$. Set
 \[ {\mathcal I}((K, \phi) X, U_{0}) = \{ f  \in IND_{H}^{G}(\tau)^{(K , \phi)} \ | \ f(X) \subseteq U_{0}, \ {\rm supp}(f) \subseteq   H \cdot X \cdot K \} , \]
 which is clearly finite-dimensional. 
 
 Now suppose that $(\tau, U)$ is ${\mathcal M}_{c}(H)$-admissible. Let $(K, \phi) \in {\mathcal M}_{c}(G)$
 and choose a finite set $X$ such that $H \cdot X \cdot K = G$. Let $L = \bigcap_{x \in X} \ xKx^{-1} $ so that $L \bigcap H$ is a compact open subspace of $H$. For each $x \in X$ we have a character $(x^{-1})^{*}(\phi) : xKx^{-1} \bigcap H \longrightarrow k^{*}$
 Take $U_{0} = {\rm Span}_{x \in X} \ U^{(L \bigcap H, (x^{-1})^{*}(\phi))}$ and suppose that $f \in IND_{H}^{G}(\tau)^{(K, \phi)}$. Therefore, if $h \in H, k \in K, x \in X$, then $f(hxk) = \phi(k) \tau(h)(f(x)$. Now suppose that we have any element of $z \in L \bigcap H$ then there is $x \in X$ and $k \in K$ such that $z = xkx^{-1} \in H$.
 Therefore $\tau(xkx^{-1})( f(x)) = f(xkx^{-1}x) = \phi(k) f(x)$ so that $f(x) \in U_{0}$. Therefore 
 $IND_{H}^{G}(\tau)^{(K, \phi)} \subseteq  {\mathcal I}((K, \phi) X, U_{0}) $ is finite-dimensional, as required.
 
 For part (v) define the map in the reverse direction
 \[ \Phi :  {\rm Hom}_{H}(V, U) \longrightarrow    {\rm Hom}_{G}(V, IND_{H}^{G}(\tau)) \]
 by $(\Phi(f)(v))(g) = f(\pi(g)v)$ which is a function of $g \in G$ which on $gk$ satisfies
 $(\Phi(f)(v))(gk) = f(\pi(g)\pi(k)v) = \phi(k)f(\pi(g)\pi(k)v) = \phi(k) (\Phi(f)(v))(g) $ if $v \in V^{(K, \phi)}$ so that the ${\mathcal M}_{c}(G)$-smoothness of $V$ ensures tha $\Phi(f)$ maps $V$ into $IND_{H}^{G}(\tau)$.
 Also composing with $\Lambda$ gives $v \mapsto f(v)$ so that $\Phi \cdot (\Lambda \cdot -) = 1$.
 Also $ (\Lambda \cdot -) \cdot \Phi  = 1$ as in the classical case (\cite{WCnotes} Theorem 2.4.1(e)  p.27).
 
 \section{Analogue of Jacquet's Theorem for ${\mathcal M}_{c}(G)$-admissibility}
 
 In this section I shall establish for ${\mathcal M}_{c}(G)$-admissibility results which are analogous to those of 
( \cite{AJS79} Chapter III, \S2.3), which is the source of any notation that I forget to elucidate.

 Let $(P,A)$ be a $p$-pair - $P = MN$. Hence $P$ is a parabolic subgroup and $N$ is its unipotent radical and $P=MN$ is its Levi decomposition.
 
For $\chi : N \longrightarrow k^{*}$ let $V_{\chi}(N)$ be the subspace of $V$ spanned by elements of the form $\pi(n)(v) - \chi(n) v$ as $v \in V$ varies\footnote{When $\chi =1$ \cite{AJS79} writes $V(N)$ sometimes as $V(P)$ but still meaning the subspace of $V$ generated by elements of the form $n \cdot v - v$ with $n \in N, v\in V$. See (\cite{AJS79} p.82).}. The proof of (\cite{AJS79} Lemma 2.2.1 p.82; also a particular case is \cite{DB96} p.461 Proposition 4.4.3)
shows that 
\begin{lemma}{$_{}$}
\label{2.1}
\begin{em}
 $v \in V_{\chi}(N)$ if and only if there is a compact open subset $N(v)$ of $N$ such that
 \[   \int_{N(v)}  \chi^{-1}(n)  \pi(n)(v) dn = 0. \]
\end{em}
\end{lemma} 

The element $m \in M$ maps $V_{\chi}(N)$ to $V_{(m^{-1})^{*}(\chi)}(N)$ where $(m^{-1})^{*}(\chi)(n) = \chi(m^{-1}nm)$ and $n \in N$ maps $V_{\chi}(N)$ to itself. In the classical case $\chi = 1$ and in that case I shall write $V(N)$ rather than $V_{1}(N)$. 

Let $\tau$ be an ${\mathcal M}_{c}(P)$-smooth representation of $P$ which is trivial on $N$ in which $P$ acts on a vector space $U$. Set $\rho = {\rm IND}_{P}^{G}(\tau)$, as in Theorem \ref{1.13}.
Therefore, by Theorem \ref{1.13}(i) $\rho$ is ${\mathcal M}_{c}(G)$-smooth and, since $P \backslash G$ is compact, by Theorem \ref{1.13}(iv) $\rho$ is ${\mathcal M}_{c}(G)$-admissible if $\tau$ is 
${\mathcal M}_{c}(P)$-admissible.

Following the method of (\cite{AJS79} Chapter III, \S2.3) we shall prove a series of lemmas in order to establish the following result :
 
 \begin{theorem}{(Analogue of Jacquet's Theorem)}
\label{2.2}
\begin{em}

Let $X_{0}$ be a non-zero ${\mathcal M}_{c}(G)$-admissible $G$-subrepresentation of ${\rm IND}_{P}^{G}(\tau)$. Then $\Delta(X_{0})$ is a non-zero  ${\mathcal M}_{c}(M)$-admissible $M$-subrepresentation of $U$.
\end{em}
\end{theorem}

Let $(P_{0}, A_{0}) $ is a minimal $p$-pair such that $(P,A) > (P_{0}, A_{0})$ and $\{ K_{j}, 1 \leq j \leq \infty \}$ is a sequence of compact open subsets of $G$ such that

(i)  $\{ K_{j}, 1 \leq j \leq \infty \}$ is a fundamental sequences of neighbourhoods of the identity in $G$,

(ii) Let $(P,A)$ be a $p$-pair of $G$ with Levi decomposition $P=MN$. If $(P,A)$ is standard with respect to $P_{0},A_{0})$ then
$K_{j} = \overline{N}_{j}M_{j}N_{j} = N_{j}M_{j}\overline{N}_{j}$ with $N_{j} = N \bigcap K_{j}, M_{j} = M \bigcap K_{j}, \overline{N}_{j} = \overline{N} \bigcap K_{j}$.

Write $U_{0} = \Delta(X_{0})$ and define, for $1 \leq j \leq \infty$, 
\[ \begin{array}{ll}
X_{0}(j, \phi) &  =  X_{0}^{(K_{j}, \phi) }  \\
\\
& = \{ f \in X_{0} \ | \  \rho(k)f = f(- \cdot k) = \phi(k) \cdot f \ {\rm for \ all } \ k \in K_{j} \}  \\
\\
& =  \{ f \in X_{0} \ | \  \phi(k)^{-1} \rho(k)f =  f \ {\rm for \ all } \ k \in K_{j} \}  .
\end{array}   \] 
Also 
\[  \begin{array}{ll} 
U_{0}(j, \phi) & = U_{0}^{(M_{j}, \phi)}   \\
\\
& = \{ u \in U_{0} \ | \  \tau(m) u  = \phi(m)u \ {\rm for \ all } \ m \in M_{j}  \} \\
\\
& =  \{ u \in U_{0} \ | \  \phi(m)^{-1} \tau(m)u  =  u \ {\rm for \ all } \ m \in M_{j}  \} 
\end{array} \]

\begin{lemma}{$_{}$}
\label{2.3}
\begin{em}
$\Delta(X_{0}(j, \phi)) \subseteq U_{0}(j, \phi)$.
\end{em} 
\end{lemma} 

{\bf Proof:}

If $m \in M_{j}$ then $m \in K_{j}$. Now, by (\cite{AJS79} Lemma 2.3.2),  $\tau(m)(f(1)) = f(m) = \phi(m) \Delta(f)$. $\Box$

\begin{lemma}{$_{}$}
\label{2.4}
\begin{em}

For any $j>0$, $\Delta(X_{0}(j, \phi))$ is stable under $\tau$ restricted to $A$.
\end{em} 
\end{lemma}

{\bf Proof:}

Since $\Delta(X_{0}(j, \phi))$ is finite-dimensional the set 
\[ S = \{a \in A \ | \ \tau(a)\Delta(X_{0}(j, \phi)) \subseteq \Delta(X_{0}(j, \phi))  \}  \]
is a group. Since $A$ is the group generated by $A^{+}(t)$ for any $t>0$, in the notation of (\cite{AJS79} Chapter 0),  it is sufficient to show that there exists $t>0$ such that $a^{-1} \in S$ for all $a \in A^{+}(t)$.

Choose $t$ such that $a \overline{N}_{j}a^{-1} \subseteq \overline{N}_{j}$ ($\overline{N}$ is the opposite of $N$, which is therefore also normalised by $M$) provided that $a \in A^{+}(t)$.

For $a \in A^{+}(t)$ and $f \in X_{0}(j, \phi)$ set
\[ \begin{array}{ll}
f_{(a,j)}  &   =  \int_{K_{j}} \phi(k)^{-1} \rho(ka^{-1}) dk f  \\
\\
 & = \int_{N_{j}} \int_{M_{j}}  \int_{ \overline{N}_{j}   }  \ (\phi(nm  \overline{n})^{-1} \rho(nm \overline{n} a^{-1})\cdot f ) dn dm d\overline{n}  ,
\end{array}  \]
which, assuming normalised measures, $f_{(a,j)} \in X(j, \phi)$ since we have averaged over $K_{j}$.

Notice that $\phi(n) = 1$ for $n \in N_{j}$, if $\Delta(f) \not= 0$, since $\tau(n)$ acts trivially on $U$.

Applying $\Delta$ we obtain, as $a$ centralises with $M_{j}$ (\cite{WCnotes}  p.11),
\[ \begin{array}{ll}
\Delta(f_{a,j}) &  =  f_{a,j}(1) \\
\\
& = \int_{M_{j} }  \int_{\overline{N}_{j} }   \phi(a)  \phi(m \overline{n} a^{-1} )^{-1} \rho(m \overline{n}a^{-1})  dm d\overline{n} f(1) \\
\\
& = \int_{ \overline{N}_{j} }  \phi(a)  d\overline{n}  f(a^{-1} (a\overline{n}a^{-1})) \\
\\
& =  \phi(a) f(a^{-1}) \ {\rm since }  \    a \overline{N}_{j} a^{-1} \subseteq \overline{N}_{j}   \\
\\
& = \phi(a) \tau(a^{-1}) \Delta(f) ,
\end{array} \]
as required. $\Box$

\begin{lemma}{$_{}$}
\label{2.5}
\begin{em}
 $U_{0}(j, \phi) \subseteq \Delta(X_{0}(j, \phi))$.
\end{em}
\end{lemma}

{\bf Proof:}

We must show that, if $u_{0} \in U_{0}(j, \phi)$, then there exists $f_{0} \in X_{0}(j, \phi)$ such that 
$\Delta(f_{0}) = u_{0}$. 

There exists $j'$ and $\phi_{i} : K_{j'} \longrightarrow k^{*}$ for $1 \leq i \leq N$, with $f_{i} \in X(j', \phi_{i})$ such that $u_{0} = \sum_{i=1}^{N} \ \Delta(f_{i})$.

Suppose that $N$ is minimal integer occurring  in the above type of relation. If $K_{j'} \subseteq K_{j}$ then $M_{j'} \subseteq M_{j}$ and for $m \in M_{j'}$ we have
\[ \begin{array}{ll}
\sum_{i=1}^{N}  \   \phi(m) f_{i}(1) &   =   \phi(m) u_{0} \\
\\
& = \tau(m) u_{0}  \\
\\
 & = \sum_{i=1}^{N} \  ( \rho(m)f_{i})(1)  \\
\\
&  =  \sum_{i=1}^{N} \ f_{i}(m) \\
\\
&  =   \sum_{i=1}^{N} \  \phi_{i}(m) f_{i}(1) . \\
\end{array} \]
Therefore we can eliminate one of the $f_{i}$'s unless each $f_{i} \in X_{0}(j', \phi)$, but such an eliniation would contradict minimality.

If $K_{j} \subset K_{j'}$ a similar minimality argument shows that we may assume $f_{i} \in X_{0}(j', \phi)$ for each $i$, which implies the result in this case.

We continue, assuming at $K_{j'} \subseteq K_{j}$ and that $f_{i} \in X_{0}(j', \phi)$ for each $i$.

Let $\lambda  \in C_{c}(M//M_{j}) $ be the characteristic function of $M_{j}$ multiplied by the reciprocal of the Haar measure of $M_{j}$ (\cite{AJS79} p.30). We have 
\[ \begin{array}{ll}
u_{0} & = \int_{M_{j}} \lambda(m) \phi(m)^{-1} \tau(m) dm u_{0} \\
\\
& =  \sum_{i=1}^{N} \  \Delta(\int_{M_{j}} \lambda(m) \phi(m)^{-1}  \rho(m) dm  f_{i}) \\
\\
&  =   \sum_{i=1}^{N} \  \int_{M_{j}} \lambda(m)  \phi(m)^{-1}  f_{i}(m) dm  . 
 \end{array}  \] 
 
Assume that $f_{i} \in X_{0}(j', \phi)$ and choose $t > 0$ so that, if $a \in A^{+}(t)$, then 
$a \overline{N}_{j} a^{-1} \subseteq \overline{N}_{j'}$. Fix $a \in A^{+}(t)$. Set $\lambda_{a}(m) = \lambda(am)$  for $m \in M$ and $\lambda'_{a}(x) = \lambda_{a}(m)$ if $x = km \in K_{j}M$, $0$ if $x \not\in K_{j}M$. Then $\lambda'_{a} \in C_{c}(K_{j} \backslash G)$, so define $f_{1,i}$ to be the convolution of $\lambda'_{a} $ with
$\phi^{-1} f_{i}$ so that  $f_{1,i} = (\lambda'_{a}) *  \phi^{-1} f_{i}  \in X_{0}(j, \phi)$. 

Then
\[ \begin{array}{ll}
((\lambda'_{a}) *  \phi^{-1} f_{i})(1) & = \int_{G}  (\lambda'_{a}(x) \phi(x)^{-1} f(x) dx \\
\\
& = \int_{K_{j}} \phi(xa^{-1})^{-1} f_{i}(xa^{-1} ) dx \\
\\
&  =  \phi(a) \int_{\overline{N}_{j}} \phi(axa^{-1})^{-1} f_{i}(axa^{-1} ) dx \\
\\
& =  \phi(a) \tau(a)  \int_{\overline{N}_{j}}  f_{i}(a \overline{n}a^{-1}) d \overline{n} \\
\\
& = \phi(a) \tau(a)(f_{i}(1))  ,
\end{array} \]
since $a \overline{N}_{j} a^{-1} \subseteq   \overline{N}_{j'}.$
$\Box$

Lemmas \ref{2.3}, \ref{2.4} and \ref{2.5} establish Theorem \ref{2.2}. In particular, combining Theorem \ref{2.2} and 
Theorem \ref{1.13}(iv) we obtain the following result.
\begin{corollary}{$_{}$}
\label{2.6} 
\begin{em}

 In the notation of \S\S\ref{2.2}-\ref{2.5}, $(\tau, U)$ is ${\mathcal M}_{c}(M)$-admissible if and only if  ${\rm IND}_{H}^{G}(\tau)$ is ${\mathcal M}_{c}(G)$-admissible.
\end{em}
\end{corollary} 
\begin{theorem}{$_{}$}
\label{2.7}
\begin{em}
 If $V$ is an ${\mathcal M}_{c}(G)$-admissible representation then $V/V(P) $ is an ${\mathcal M}_{c}(M)$-admissible representation\footnote{Recall that, following the ambiguous notation of \cite{AJS79}, $V(P)$ is also denoted by $V(N)$.}.
\end{em}
\end{theorem}

{\bf Proof:}

By Frobenius reciprocity (see Theorem \ref{1.13}(v)) we have an isomorphism of the form
  \[  ( \Lambda \cdot - ) : {\rm Hom}_{G}(V, IND_{P}^{G}(V/V(P) )) \stackrel{\cong}{\longrightarrow }  {\rm Hom}_{P}(V, V/V(P) ) \]
 and we denote by $T$ the $G$-map on the left which corresponds the the canonical quotient of $P$-representations on the right. Setting $X_{0} = T(V) \subset  IND_{P}^{G}(V/V(P) )$, which is 
${\mathcal M}_{c}(G)$-admissible by Corollary \ref{1.9},  implies that $\Delta(X_{0}) = V/V(P)$ is an ${\mathcal M}_{c}(M)$-admissible representation, by Theorem \ref{2.2}. $\Box$

\section{Related remarks and questions}
If $\chi$ is a character of $P$ which is trivial on $M$ then it gives an $M$-fixed character on $N$ and conversely an an $M$-fixed character on $N$ $\chi$ extends to $\chi(mn) = \chi(n)$) on $P$, which is trivial on $M$.

In this case $V/V_{\hat{\chi}}(P)$ of Lemma \ref{2.1} is a $P$-representation\footnote{Similarly, recall that, following \cite{AJS79}, we adopt the ambiguous notation of writing $V_{\hat{\chi}}(P)$ and $V_{\hat{\chi}}(N)$
for the same subspace of $V$.}.

If $\chi: P \longrightarrow k^{*}$ is a continuous character which is trivial on $M$ and $(U, \tau)$ is an ${\mathcal M}_{c}(P)$-admissible representation then so is $U' = U \otimes \chi$. Therefore, by Theorem \ref{1.13}(iv), $X' = {\rm IND}_{P}^{G}(U \otimes \chi)$ is  an ${\mathcal M}_{c}(G)$-admissible representation.

In \cite{Sn18}, \cite{Sn20}, \cite{Sn20b} and \cite{Sn20c} one is often concerned with the subcategory of admissible representations with a fixed central character $\underline{\phi}$ and ${\mathcal M}_{cmc. \underline{\phi}}(G)$-admissibility, where  ${\mathcal M}_{cmc. \underline{\phi}}(G)$ is the poset of pairs $(H, \phi)$ where $Z(G) \subseteq H$ and $\phi$ is a continuous character extending $\underline{\phi}$ and $H$ is compact open modulo the centre. There is an analogue of $  IND_{P}^{G}(-)$ in this setting and presumably the obvious analogue of the classical Jacquet theory.

For the classical definition can one modify the Jacquet proof to show $V$ $G$-admissible implies $V/V_{\chi}(P)$ is $P$-admissible?

If that works then can one try to prove it for ${\mathcal M}_{c}(G)$-admissibility?

A question left in suspense in \cite{Sn20} was the derivation of the general formula for generators of the hyperHecke algebra in term of convolution products. I gave the formula in the case when all characters on compact opens have finite image (as in the classical case). It seems to me that, with a little more expertise with Mahler's theorem (see \cite{AMR00}) my previous proof might be promoted to the general case? Unfortunately, (see  \cite{JMOS18} and \cite{XX}) this may take rather too long, in the absence of access to 
``big-hitting Langlands professionals''. This one will have to wait until the unlikely event of my next trip to Toronto, Paris or Princeton.

There is another dangling, incompletely sketched problem, vaguely conjectured in \cite{Sn20b}, concerning the Hopf-like algebra structure of the hyperHecke algebra of general linear groups. An additional discussion of this algebra is skectched in \cite{Sn20c}.


\begin{thebibliography}{PD84}

\bibitem{BZ76}  J. Bernstein and A. Zelevinski: Representations of the group $GL(n,F)$ where $F$ is a local non-Archimedean field; Uspekhi Mat. Nauk. {\bf 31} 3 (1976) 5-70.

\bibitem{BZ77}  J. Bernstein and A. Zelevinski: Induced representations of reductive $p$-adic groups I; Ann. ENS {\bf 10} (1977) 441-472.

\bibitem{ABJT65} A. Borel and J. Tits: Groupes r\'{e}ductifs; I.H.E.S. Pub. Math.  no. 27 (1965) 55-150 MR {\bf 34} \#7527.

\bibitem{FBr56} F. Bruhat: Sur les repr\'{e}sentations induites des groupes de Lie; Bull. Soc. Math. France 84 (1956) 97-205.

\bibitem{FBr61}  F. Bruhat: Distributions sur un groupe localement compact et applications \`{a} l'\'{e}tude des repr\'{e}sentations des groupes $p$-adiques; ; Bull. Soc. Math. France 89 (1961) 43-75.

  
\bibitem{DB96} Daniel Bump: {\em Automorphic forms and representations}; Cambridge studies in advanced math. {\bf 55} (1998).

\bibitem{BH06}  C.J. Bushnell and G. Henniart: {\em The Local Langlands Conjecture for $GL(2)$}; Grund. Math. Wiss. \#335; Springer Verlag (2006).

\bibitem{WCnotes} W. Casselman: Introduction to  the theory of admissible representations of $p$-adic groups; Preprint.
  
\bibitem{PD84}  P. Deligne: Le ``centre'' de Bernstein;  {\em Repr\'{e}sentations des groupes r\'{e}ductifs sur un corps local}  Travaux en cours, Hermann, Paris (1984) 1-32.
 
 \bibitem{JAG55} J.A. Green: The characters of the finite general linear groups; Trans. Amer. Math. Soc. 80 (1955) 402-447.
 
  \bibitem{HC70} Harish-Chandra: Harmonic analysis on reductive $p$-adic groups; Lecture Notes in Math. \#162, (notes by G. van Dijk) Springer Verlag (1970).
 
 \bibitem{HC74} Harish-Chandra: Harmonic analysis on reductive $p$-adic groups; Proc. Symp. Pure Math. {\bf 26} (ed. C.C. Moore) Amer. Math. Soc. (1974) 167-192.
 
 \bibitem{JMOS18} J.M. O'Sullivan and C.N. Harrison: Myelofibrosis: Clinicopathologic Features, Prognosis and Management; Clinical Advances in Haematology and Oncology {\bf 16} (2) February 2018.
 
 \bibitem{XX} Clinical articles: https://www.ncbi.nlm.nih.gov/pmc/articles/PMC6514804, https://www.mayoclinicproceedings.org/article/S0025-6196(17)30380-4/pdf and https://www.hematologyandoncology.net/archives/july-2019/management-of-advanced-phase-myeloproliferative-neoplasms.

 
 \bibitem{RPL} R.P. Langlands: On the classification of irreducible representations of real algebraic groups; preprint.
  
\bibitem{IGM80}  I.G. Macdonald: Zeta functions attached to finite general linear groups; Math. Annalen 249 (1980) 1-15.

\bibitem{MZ55} D. Montgomery and L. Zippin: {\em Topological Transformation Groups}; Interscience New York (1955).

\bibitem{AMR00} A. M. Robert: {\em A Course in $p$-adic Analysis}; Grad. Texts in Math. $\#$198, Springer Verlag (2000). 

\bibitem{AJS78} Allan J. Silberger: The Langlands quotient theorem for p-adic groups; Math. Annalen 236 no. 2 (1978) 95-104.

\bibitem{AJS79}  Allan J. Silberger: {\em Introduction to harmonic analysis on $p$-adic reductive groups}; Math. Notes Princeton Unviersity Press (1979).

\bibitem{Sn94}  V.P. Snaith:  {\em Explicit Brauer Induction (with applications to algebra and number theory)}; Cambridge studies in advanced mathematics \#40,Cambridge University Press (1994).
  
\bibitem{Sn18} V.P. Snaith:  {\em Derived Langlands}; World Scientific (2018).

  
  \bibitem{Sn20} V.P. Snaith: Derived Langands II; arXiv: 3100675 [math.RT] 24 Mar 2020.
  
\bibitem{Sn20b}   V.P. Snaith:  Derived Langlands III: PSH algebras and their numerical invariants; ArXiv June 2020.
  
  \bibitem{Sn20c}  V.P. Snaith: The Hopflike properties of the hyperHecke algebra; 2020 preprint University of Sheffield homepage (preliminary version).

  
 \end{thebibliography}
 \end{document}